\documentclass[11pt]{article}
\usepackage{etex}
\usepackage[a4paper]{geometry}

\usepackage{amsmath}
\usepackage{amssymb}
\usepackage{amsthm}
\usepackage{mathrsfs}
\usepackage{bbm}
\usepackage{empheq}
\usepackage[loose]{subfigure}
\usepackage{epsfig}
\usepackage{graphicx}
\usepackage{psfrag}
\usepackage[usenames,dvipsnames]{pstricks}
\usepackage{pst-plot}
\usepackage[colorlinks]{hyperref}
\usepackage{tabls}
\usepackage{paralist}
\usepackage{hyperref}
\DeclareSymbolFontAlphabet{\Bbb}{AMSb}
\usepackage[all]{xy}
\newlength{\fixboxwidth}
\newcommand{\COMMENT}[1]{}
\newcommand{\E}{\mathbb{E}}

\newcommand{\R}{\mathbb{R}}

\newcommand{\quark}{\setbox0\hbox{$x$}\hbox to\wd0{\hss$\cdot$\hss}}

\newcommand{\overbar}[1]{\mkern 1.5mu\overline{\mkern-1.5mu#1\mkern-1.5mu}\mkern 1.5mu}

\newcommand{\s}{\sigma}

\newcommand{\snorm}[1] {\Vert #1 \Vert}

\newtheorem{thm}{Theorem}[section]

\newtheorem{cor}[thm]{Corollary}
\theoremstyle{definition}

\newtheorem{rmk}[thm]{Remark}

\hypersetup{
    linkcolor=blue,
}

\title{Conditioning Gaussian measure on Hilbert space}
\author{Houman Owhadi and Clint Scovel
\\
California Institute of Technology
}
\date{\today}

\makeatletter
\@addtoreset{equation}{section}
\@addtoreset{figure}{section}
\@addtoreset{table}{section}
\makeatother

\renewcommand{\thefigure}{\arabic{section}.\arabic{figure}}

\makeatletter
\renewcommand{\p@subfigure}{\thefigure}
\makeatother

\newcounter{mycount}

\begin{document}

\maketitle
\abstract{
For a Gaussian measure on a separable Hilbert space with covariance
operator $C$, we show  that the family of conditional measures associated with conditioning
on a closed subspace  $S^{\perp}$ are Gaussian with covariance operator the short $\mathcal{S}(C)$ of the operator $C$ to $S$. We provide two proofs. The first uses the theory of Gaussian Hilbert spaces and
a characterization of the shorted operator by Andersen and Trapp.  The second uses
 recent developments by Corach, Maestripieri and Stojanoff
on the relationship between the shorted operator and $C$-symmetric oblique projections onto $S^{\perp}$.
 To obtain the assertion when such
projections do not exist, we develop an approximation result for the shorted operator by showing,
  for any positive operator $A$, how to construct a sequence of approximating
operators
 $A^{n}$  which possess  $A^{n}$-symmetric oblique projections onto $S^{\perp}$ such that
the sequence of shorted operators $\mathcal{S}(A^{n})$ converges to $\mathcal{S}(A)$ in the weak operator
 topology.  This result combined with the martingale convergence of  random variables associated
with the corresponding approximations $C^{n}$ establishes the main assertion in general.
 Moreover, it in turn strengthens
the approximation theorem for shorted operator when the operator is trace class; then the sequence
of shorted operators $\mathcal{S}(A^{n})$ converges to $\mathcal{S}(A)$ in trace norm.
}
	


\section{Introduction}
\label{sec_intro}
For a Gaussian measure $\mu$  with injective covariance operator $C$ on a direct sum of finite dimensional Hilbert spaces $H=H_{1}\oplus H_{2}$, the conditional
measure associated with conditioning on the value of $H_{2}$ can be computed
in terms of the Schur complement corresponding to the partitioning of the covariance matrix $C$, see Cottle
\cite{cottle1974manifestations} for a review.  Evidently, the natural extension to infinite dimensions
of the  Schur complement is the {\em shorted operator}, first discovered by  Krein \cite{krein1947theory}
 and
 developed in Anderson and Trapp  \cite{anderson1975shorted} based on results on operator ranges
of Douglas \cite{douglas1966majorization} and Fillmore and Williams \cite{fillmore1971operator}.
   For related results, and a history,  see Pekarev \cite{pekarev1992shorts}.
  However, the connection between the shorted operator and the covariance operator
of the conditional Gaussian measure on an infinite dimensional Hilbert space appears yet to be established.
 Indeed,  Hairer, Stuart, Voss, and Wiber
\cite[Lem.~4.3]{hairer2005analysis}, see also   Stuart \cite[Thm.~6.20]{Stuart:2010},
 characterizes the conditional measure through
 a measurable extension
result of Dalecky and Fomin
\cite[Thm.~II.3.3]{dalet︠s︡kiui1991measures} of an operator defined on the Cameron-Martin reproducing kernel Hilbert space.  For other representations, see  Mandelbaum \cite{mandelbaum1984linear}, LaGatta \cite{lagatta2013continuous}, and
Tarieladze and Vakhania's \cite{tarieladze2007disintegration} extension of the  optimal linear approximation
 results of
Lee and Wasilkowski \cite{lee1986approximation} from finite to infinite rank.
Tarieladze \cite{tarieladzeinformation} asserts that
this latter result extends one in
the Information-Based Complexity of
 Traub, Wasilkowski and Wozniakowski \cite{traubniakowski} which is relevant to
Grid Computing  as described in
 Foster and Kesselman \cite{foster2003grid}.
 The primary  purpose of this paper is to instead represent the conditional measure in terms of the shorted operator. We provide two distinct proofs of this representation. The first uses the theory of Gaussian Hilbert spaces 
and
a characterization of the shorted operator by Andersen and Trapp. The second proof, corresponding 
to the secondary 
purpose of this paper,
 uses recent developments by Corach, Maestripieri and Stojanoff
on the relationship between the shorted operator and $A$-symmetric oblique projections. This latter approach has the advantage that it
facilitates a general  approximation technique that not only can be used to approximate the covariance operator but the conditional expectation operator. This is accomplished through the development of an approximation theory for the shorted operator in terms of oblique projections followed by an application  of the martingale convergence
theorem. 
Although the proofs are not fundamentally difficult, the result (which appears to have been missed in the literature) provides a simple characterization of the conditional measure, leading to significant approximation results.
For instance, the attainment of the main result through the martingale approach feeds back a strengthening of the approximation  theorem  for the shorted operator that was developed for that purpose:  when the operator is trace class the approximation improves from weak convergence to convergence in trace norm.

Let us review the basic results on Gaussian measures on Hilbert space.
 A measure $\mu$ on a Hilbert space $H$ is said to be Gaussian if, for each $h \in H$ considered as a continuous linear function
$h:H \rightarrow \R$ by $h(x):=\langle h,x\rangle, x \in H$, we have that the pushforward measure
$h_{*}\mu$ is  Gaussian,
 where we say that a Dirac measure is Gaussian.
For a Gaussian measure $\mu$, its mean $m$  is defined by
\[ \langle h ,m\rangle :=\int_{H}{\langle h, x\rangle d\mu(x)},  \quad  h \in H \, \]
and its covariance operator $C:H \rightarrow  H $ is defined by
\[ \langle Ch_{1},h_{2}\rangle :=\int_{H}{\langle h_{1}, x\rangle \langle h_{2}, x\rangle d\mu(x)}-
\langle  h_{1}, m\rangle \langle  h_{2}, m\rangle ,  \quad h_{1}, h_{2} \in H \, .\]
A Gaussian measure has a well defined mean
and a continuous covariance operator, see e.g.~Da Prato and Zabczyk \cite[Lem.~2.14]{da2014stochastic}.
Finally,
 Mourier's Theorem \cite{mourier1953elements},
 see Vakhania, Tarieladze and Chobanyan \cite[Thm.~IV.2.4]{vakhania1987probability}, asserts,
 for any $m \in H$ and any positive symmetric trace class operator $C$,  that there exists a Gaussian measure with mean $m$
and covariance operator $C$, and that all Gaussian measures have a well defined mean and  positive trace class covariance operator.  This characterization also follows from Sazonov's Theorem \cite[Thm.~1]{sazonov1958remark}.

Since separable Hilbert spaces are Polish, it follows from
the product space version, see e.g.~Dudley \cite[Thm.~10.2.2]{Dudley:2002},
of the theorem on the existence and uniqueness of regular conditional probabilities
on Polish spaces, that any Gaussian measure $\mu$
on a direct sum $H=H_{1}\oplus H_{2}$ of separable Hilbert spaces has  a regular conditional probability, that is
there is a family $\mu_{t}, t \in H_{2}$ of conditional  measures corresponding to conditioning on $H_{2}$.
 Moreover,
 Tarieladze and Vakhania \cite[Thm.~3.11]{tarieladze2007disintegration}  demonstrate
that the corresponding family of  conditional measures are Gaussian.
Bogachev's  \cite[Thm.~3.10.1]{bogachev1998gaussian}  theorem of normal correlation of Hilbert space valued
Gaussian random variables shows that if two Gaussian random vectors $\xi$ and $\eta$ on a separable Hilbert space
$H$ are jointly Gaussian in the product space, then
$\E[\xi|\eta]$ is a Gaussian random vector and
$\xi=\E[\xi|\eta]+\zeta$ where $\zeta$ is Gaussian random vector which is independent of $\eta$. Consequently,
for any two vectors $h_{1}, h_{2} \in H$ we have
\begin{eqnarray*}
\E\Bigl[\bigl\langle \xi-\E[\xi|\eta], h_{1}\bigr\rangle\bigl\langle \xi-\E[\xi|\eta], h_{2}\bigr\rangle\big|
\eta\Bigr]&=&\E\Bigl[\bigl\langle\zeta , h_{1}\bigr\rangle\bigl\langle \zeta, h_{2}\bigr\rangle\big|
\eta\Bigr]\\
&=&\E\Bigl[\bigl\langle\zeta, h_{1}\bigr\rangle\bigl\langle \zeta, h_{2}\bigr\rangle\Bigr]
\end{eqnarray*}
and so we  conclude that, just as in the finite dimensional case,
 the conditional covariance operators are independent of the values of the conditioning
 variables.

Since both proof techniques will utilize
the characterization of conditional expectation as orthogonal projection, we introduce these notions now.
Consider the Lebesgue-Bochner space $L^{2}(H,\mu,\mathcal{B}(H))$  space of (equivalence classes) of $H$-valued
Borel measurable functions  on $H$ whose squared norm
\[ \snorm{f}^{2}_{L^{2}(H,\mu,\mathcal{B}(H))}:=\int_{H}{\snorm{f(x)}^{2}_{H}d\mu(x)}\ \]
is integrable.
For a sub $\s$-algebra $\Sigma \subset \mathcal{B}(H)$ of the Borel $\s$-algebra, consider the
corresponding
Lebesgue-Bochner space $L^{2}(H,\mu,\Sigma)$.
As in the scalar case, one can show that
 $L^{2}(H,\mu,\mathcal{B}(H))$  and $L^{2}(H,\mu,\Sigma)$ are Hilbert spaces  and that
$L^{2}(H,\mu,\Sigma) \subset L^{2}(H,\mu,\mathcal{B}(H))$ is a closed subspace.
Then, if we note
that  contractive projections on Hilbert space are orthogonal, see
Rao  \cite[Rmk.~9, pg.~51]{rao2011foundations}, it follows from
Sundaresan \cite[Prop.~4]{Sundaresan}, see Diestel and Uhl \cite[Thm.~V.1.4]{DiestelUhl},
that conditional expectation amounts to orthogonal  projection.

\section{Shorted Operators}
A symmetric operator $A:H \rightarrow H$ is called {\em positive} if
$\langle Ax, x\rangle \geq 0$ for all $x \in H$. We denote by $L_{+}(H)$ the set of positive  operators  and we denote  such positivity by
$A \geqq 0$. Positivity induces the (L\"{o}wner) partial order $\geqq $ on $L_{+}(H)$.
 For a closed subspace $S \subset H$ and a positive operator $A\in L_{+}(H)$  consider the set
\[ \mathcal{H}(A,S):=\bigl\{X \in L_{+}(H):  X \leqq A\,\,  \text{and} \,\,  R(X)\subset S \bigr\}\, .\]
Then, according to Pekarev \cite{pekarev1992shorts}, Krein \cite{krein1947theory} and later Anderson and Trapp  \cite{anderson1975shorted} showed that
$\mathcal{H}(A,S)$ contains a maximal element, which we denote by $\mathcal{S}(A)$ and call the {\em short}
of $A$ to $S$. For another closed subspace $T\subset H$, we denote the short of $A$ to $T$ by
$\mathcal{T}(A)$.
  In the proof,
 Anderson and Trapp  \cite{anderson1975shorted} demonstrate that when $A$ is invertible, that in terms of its
 $(S,S^{\perp})$
partition representation
\[A= \begin{pmatrix}
  A_{SS} &  A_{SS^{\perp}}  \\
  A_{S^{\perp}S} & A_{S^{\perp}S^{\perp}}
 \end{pmatrix}
\]
that $A_{S^{\perp}S^{\perp}}$ is invertible and
\[ \mathcal{S}(A)= \begin{pmatrix}
A_{SS}-A_{SS^{\perp}}A_{S^{\perp}S^{\perp}}^{-1} A_{S^{\perp}S}  &  0  \\
  0 & 0
 \end{pmatrix}\, .
\]
It is easy to show that the assertion holds under the weaker assumption
  that $A_{S^{\perp}S^{\perp}}$ be invertible.
Moreover,  Anderson and Trapp  \cite[Cor.~1]{anderson1975shorted} asserts  for $A, B \in L_{+}(H)$, that
\[
A \leqq B \implies \mathcal{S}(A) \leqq\mathcal{S}(B) \, ,\]
 that is, $\mathcal{S}$ is monotone in the L\"{o}wner
ordering. In addition,  \cite[Cor.~5]{anderson1975shorted}  asserts that
for two closed subspaces $S$ and $T$, we have
\[
 \mathcal{(S\cap T)}(A) =\mathcal{S}\bigl(\mathcal{T}(A)\bigr)\, .\]
Finally, \cite[Thm.~6]{anderson1975shorted} asserts that
 if $A:H \rightarrow H$ is a positive operator and $S \subset H$ is a closed
 linear subspace, then
\begin{equation}
\label{id_short}
  \bigl\langle \mathcal{S}(A)s,s\bigr\rangle =\inf\Bigl\{\Bigl\langle A
\begin{pmatrix}
  s\\
t
 \end{pmatrix}, \begin{pmatrix}
  s\\
t
 \end{pmatrix} \Bigr\rangle, \, t \in S^{\perp} \Bigr\},\quad \forall s \in S\, .
\end{equation}
In Section  \ref{sec_proof1} we demonstrate how the characterization \eqref{id_short} of the shorted operator
combined with the theory of Gaussian Hilbert spaces provides a natural  proof of our main result, the following theorem. Here we consider direct sum split $H=H_{1}\oplus H_{2}$, and let
$S=H_{1}$ and $S^{\perp}=H_{2}$, so that the short $\mathcal{S}(A)$ of an operator
to the subspace $S=H_{1}$ will be written as $\mathcal{H}_{1}(A)$.
\begin{thm}
\label{thm_main0}
Consider a  Gaussian measure $\mu$ on an orthogonal direct sum $H=H_{1}\oplus H_{2}$ of separable Hilbert spaces
   with mean $m$
 and covariance operator $C$.
Then for all $t \in H_{2}$,  the conditional measure $\mu_{t}$ is a Gaussian measure
 with
 covariance operator $\mathcal{H}_{1}(C).$
\end{thm}

\section{Oblique Projections}
In this section, we will prepare for  an alternative proof of Theorem \ref{thm_main0} using oblique projections
along with the development of approximations of the covariance operator and the conditional expectation operator
 generated by natural sequences of oblique projections.
To that end, let us introduce some notations.
For a separable Hilbert space $H$,
we denote the usual, or strong, convergence of sequences by $h_{n} \rightarrow h$
and the weak convergence by $h_{n} \xrightarrow{\omega} h$.
Let $L(H)$ denote the Banach algebra of bounded linear operators on $H$.
For an operator $A \in L(H)$, we let $R(A)$ denote its range and $ker(A)$ denote its nullspace.
Recall the uniform operator topology on $L(H)$ defined by the metric
$\snorm{A}:=\sup_{\snorm{h}\leq 1}{\snorm{Ah}}.$
 We say that a sequence of operators
$A_{n} \in L(H)$ converges strongly to $A \in L(H)$, that is
\[A= s\mbox{-}\\lim_{n\rightarrow\infty}{A^{n}}\]
if $A_{n}h \rightarrow Ah$ for all $h \in H$, and we
say that $A_{n}\rightarrow A$ weakly  or
\[A= \omega\mbox{-}\\lim_{n\rightarrow\infty}{A^{n}}\]
if $A_{n}h \xrightarrow{\omega} Ah$
 for all $h \in H$.
Recall that an  operator $A \in L(H)$  is called
{\em trace class} if the {\em trace norm}
\[\snorm{A}_{1}:=\sum_{i=1}^{\infty}{\langle |A|e_{i},e_{i}\rangle }\, \]
 is finite for some orthonormal basis, where
$|A|:=\sqrt{A^{*}A}$ is the absolute value. When it is finite,  then
$tr(A):=\sum_{i=1}^{\infty}{\langle Ae_{i},e_{i}\rangle} $ is well defined, and  for all
positive trace class operators $A$ we have $tr(A)=\snorm{A}_{1}$.
The trace norm $\snorm{\cdot}_{1}$  makes the subspace $L_{1}(H) \subset L(H)$ of trace class operators
into a Banach space.
 It is well known that  the sequence of operator topologies
\[  \text{weak} \rightarrow \text{strong} \rightarrow \text{ uniform operator} \rightarrow \text{trace norm}\]
increases from left to right in strength.

 For a  positive operator $A:H \rightarrow H$, let us define the set of ($A$-symmetric) oblique projections
\[\mathcal{P}(A,S^{\perp}):=\bigl\{ Q \in L(H): Q^{2}=Q, \, R(Q)=S^{\perp},\,\, AQ=Q^{*}A\bigr \}\]
onto $S^{\perp}$, where $Q^{*}$ is the adjoint of $Q$ with respect to the scalar product $\langle  \cdot, \cdot\rangle$ on $H$. The pair $(A,S^{\perp})$ is said to be {\em compatible}, or $S^{\perp}$
is said to be {\em compatible with} $A$,  if $\mathcal{P}(A,S^{\perp})$ is nonempty.  For any oblique projection
 $Q \in \mathcal{P}(A,S^{\perp})$,
Corach, Maestripieri and Stojanoff~\cite[Prop.~4.2]{corach2001schur}
 asserts that for
  $E:=1-Q$,  we have
\begin{equation}
\label{e1}
\mathcal{S}(A)=AE=E^{*}AE\, .
\end{equation}
Moreover, when $(A,S^{\perp})$ is compatible, according to
  Corach, Maestripieri and Stojanoff \cite[Def.~3.4]{corach2001schur}, there is a special element $Q_{A,S^{\perp}} \in  \mathcal{P}(A,S^{\perp})$ defined in the following way: by \cite[Prop.~3.3]{corach2001schur}
and  the factorization theorem \cite[Thm.~2.2]{corach2001schur}
 of Douglas \cite{douglas1966majorization} and Fillmore and Williams \cite{fillmore1971operator}, there is
a unique operator
 $\hat{Q}:S \rightarrow S^{\perp}$ which satisfies
 $A_{S^{\perp}S^{\perp}}\hat{Q}=A_{S^{\perp}S}$ such that  $ker(\hat{Q})=ker(A_{S^{\perp}S})$
and $R(\hat{Q})\subset \overbar{R(A_{S^{\perp}S^{\perp}})}$.  Defining
\begin{equation}
\label{id_special}
Q_{A,S^{\perp}}= \begin{pmatrix}
  0& 0\\
\hat{Q}& 1
 \end{pmatrix}\, ,
\end{equation}
  \cite[Thm.~3.5]{corach2001schur} asserts that $Q_{A,S^{\perp}}  \in \mathcal{P}(A,S^{\perp})$.

When
 the pair $(A,S^{\perp})$ is not compatible,
 we
 seek an approximating sequence $A^{n}$ to $A$ which is compatible with
$S^{\perp}$, such that the limit of $\mathcal{S}(A^{n})$ is $\mathcal{S}(A)$.
Although  Anderson and Trapp  \cite[Cor.~2]{anderson1975shorted}  show that if $A^{n}$ is a monotone decreasing
sequence of positive operators which converge strongly to $A$, that
the decreasing sequence of positive operators $\mathcal{S}(A^{n})$ strongly converges to
$\mathcal{S}(A)$, the approximation from above by $A^{n}:=A+\frac{1}{n}I$ determines  operators
which are not trace class, so is not useful for the approximation problem for the covariance operators
for Gaussian measures. Since the trace class operators are well approximated from below by finite rank operators
one might hope to approximate $A$ by an increasing sequence of finite rank operators. However,
  it is easy to see that, in general, the same convergence result does not hold for increasing sequences.
The following theorem demonstrates, for any positive operator $A$, how to produce  a sequence of
positive operators $A^{n}$ which are compatible with $S^{\perp}$ such that
$\mathcal{S}(A^{n})$ weakly converges
to $\mathcal{S}(A)$.

 Henceforth we consider a direct sum split $H=H_{1}\oplus H_{2}$, and let
$S=H_{1}$ and $S^{\perp}=H_{2}$, so that the short $\mathcal{S}(A)$ of an operator
to the subspace $S=H_{1}$ will be written as $\mathcal{H}_{1}(A)$.
Let us also denote by $P_{i}:H\rightarrow H$ the orthogonal projections onto $H_{i}$, for $i=1,2$,
and let $\Pi_{i}:H\rightarrow H_{i}$ denote the corresponding projections
and $\Pi^{*}_{i}:H_{i} \rightarrow H$ the corresponding injections.  For any operator $A:H\rightarrow H$,
consider the decomposition
\[A= \begin{pmatrix}
  A_{11}& A_{12}\\
A_{21}& A_{22}
 \end{pmatrix}
\]
where the components are defined by
$A_{ij}:=\Pi_{i}A\Pi^{*}_{j},\,  i,j=1,2$.

\begin{thm}
\label{thm_short}
Consider a  positive operator $A:H \rightarrow H$ on a separable Hilbert space.
Then for any orthogonal split $H=H_{1}\oplus H_{2}$,  and any ordered orthonormal basis
of $H_{2}$, we let $H^{n}_{2}$ denote the span of the first $n$ basis elements and
 let $P^{n}:=P_{H_{1}}+P_{H^{n}_{2}}$ denote the orthogonal projection onto $H_{1}\oplus H^{n}_{2}$.  Then
the sequence of positive operators
\[A^{n}:=P^{n}AP^{n},\quad n=1, \ldots\]
is compatible with $H_{2}$ and
 \[\mathcal{H}_{1}(A)= \omega\mbox{-}\\lim_{n\rightarrow\infty}{ \mathcal{H}_{1}(A^{n})}\, .\]
\end{thm}
\begin{rmk}
For an increasing sequence $A^{n}$ of positive operators converging strongly to $A$, the monotonicity
of the shorting operation implies that the sequence
$\mathcal{H}_{1}(A^{n})$ is increasing, and therefore Vigier's Theorem,
 see e.g.~Halmos \cite[Prb.~120]{halmos1982hilbert}, implies that the sequence $\mathcal{H}_{1}(A^{n})$  converges
strongly.
Although
the sequence $A^{n}:=P^{n}AP^{n}$ defined in Theorem \ref{thm_short} is positive and converges
strongly to $A$, in general, it is not increasing
in the L\"{o}wner order, so that Vigier's Theorem  does not apply,  possibly suggesting why
we  only obtain convergence in the weak operator topology.  With stronger assumptions on the operator
$A$ and a well chosen selection of an ordered orthonormal basis of $H_{2}$, we conjecture
that convergence in a stronger topology may be available. In particular, as a corollary to our main
result,  when $A$ is trace class, we establish in Corollary \ref{cor_short} that
 \[ \mathcal{H}_{1}(A^{n}) \rightarrow \mathcal{H}_{1}(A) \, \, \text{in trace norm}\, . \]
\end{rmk}

For any
$m \in H$, we let $m=(m_{1},m_{2})$ denote its decomposition in $H=H_{1}\oplus H_{2}$. Moreover, for any
projection $Q:H\rightarrow H$
 with $R(Q)=H_{2}$ we let
$\hat{Q}: H_{1} \rightarrow H_{2}$ denote the unique operator such that
\[Q= \begin{pmatrix}
  0& 0\\
\hat{Q}& 1
 \end{pmatrix}\, ,
\]
and denote by $\hat{Q}^{*}:H_{2} \rightarrow H_{1}$ the adjoint of $\hat{Q}$ defined by the relation
$\langle \hat{Q}^{*}h_{2}, h_{1}  \rangle_{H_{1}} =\langle h_{2}, \hat{Q}h_{1}  \rangle_{H_{2}}$
for all $h_{1} \in H_{1}, h_{2}\in H_{2}$.

The following theorem constitutes an expansion of our main  result, Theorem \ref{thm_main0}, to include
 natural approximations for the conditional covariance operator and the conditional expectation operator.
\begin{thm}
\label{thm_main}
Consider a  Gaussian measure $\mu$ on an orthogonal direct sum $H=H_{1}\oplus H_{2}$ of separable Hilbert spaces
   with mean $m$
 and covariance operator $C$.
Then for all $t \in H_{2}$,  the conditional measure $\mu_{t}$ is a Gaussian measure
 with
 covariance operator $\mathcal{H}_{1}(C).$

If the covariance operator $C$ is compatible with $H_{2}$, then
for any oblique projection $Q$ in $\mathcal{P}(C,H_{2})\neq \emptyset$,
 the mean $m_{t}$ of the conditional measure $\mu_{t}$ is
\[m_{t}= \begin{pmatrix}
  m_{1}+\hat{Q}^{*}(t-m_{2}) \\
t
 \end{pmatrix}\, .
\]
In the general case,  for any  ordered orthonormal basis
for $H_{2}$, let $H^{n}_{2}$ denote the span of the first $n$ basis elements,
 let $P^{n}:=P_{H_{1}}+P_{H^{n}_{2}}$ denote the orthogonal projection onto $H_{1}\oplus H^{n}_{2}$,
and define the approximate  $C^{n}:=P^{n}CP^{n}$. Then
$C^{n}$ is compatible with $H_{2}$ for all $n$, and
 for any sequence $Q_{n} \in \mathcal{P}(C^{n},H_{2}) \neq \emptyset$ of
 oblique projections,
we have
\[m_{t}=  \begin{pmatrix}
  m_{1}+ \lim_{n\rightarrow\infty}{\hat{Q}_{n}^{*}(t-P_{H^{n}_{2}}m_{2})} \\
t
 \end{pmatrix}
\]
for
$\mu$-almost every $t$.
If the sequence  $Q_{n}$ eventually becomes the special element $Q_{n}=Q_{C^{n},H_{2}}$ defined
near \eqref{id_special},  then we have
\[m_{t}=  \begin{pmatrix}
  m_{1}+ \lim_{n\rightarrow\infty}{\hat{Q}_{n}^{*}(t-m_{2})} \\
t
 \end{pmatrix}
\]
for
$\mu$-almost every $t$.
\end{thm}

As a corollary to Theorem \ref{thm_main}, we  obtain a strengthening of the assertion
of Theorem \ref{thm_short} when the operator $A$ is trace class.
\begin{cor}
\label{cor_short}
Consider the situation of Theorem \ref{thm_short} with $A$ trace class. Then
\[ \mathcal{H}_{1}(A^{n}) \rightarrow \mathcal{H}_{1}(A)\,\, \text{in trace norm}\, . \]
\end{cor}

\section{Proofs}
\subsection{First proof of Theorem \ref{thm_main0}}
\label{sec_proof1}
 Consider the Lebesgue-Bochner space $L^{2}(H,\mu,\mathcal{B}(H))$  space of (equivalence classes) of $H$-valued
Borel measurable functions on $H$ whose squared norm
\[ \snorm{f}^{2}_{L^{2}(H,\mu,\mathcal{B}(H))}:=\int_{H'}{\snorm{f(x)}^{2}_{H}d\mu(x)}\ \]
is integrable.
For any square Bochner
integrable function $f \in L^{2}(H,\mu,\mathcal{B}(H))$ and any
 $h \in H$, we have that
$\langle f,h \rangle $ is square integrable, that is  $\langle f,h \rangle \in L^{2}(\R,\mu,\mathcal{B}(H))$.
 Moreover, it is easy to see, see e.g.~
\cite[Lem.~11.45]{AliprantisBorder:2006},  that if $f$ is Bochner integrable, then for all
$h \in H$, we have $\langle f,h \rangle $ is Bochner integrable and
$\int{\langle f,h \rangle d\mu}=\langle \int{fd\mu},h \rangle$.

Now consider the orthogonal decomposition $H=H_{1}\oplus H_{2}$ and
the Borel $\s$-algebra $\mathcal{B}(H_{2})$. Let us denote the shorthand notation
\[\mathcal{B}:=\mathcal{B}(H), \quad
\mathcal{B}_{2}:= \{ (H_{1},T):  T \in \mathcal{B}(H_{2})\}\, .\]
The definition of conditional expectation in  Lebesgue-Bochner space, that is that
$\E[f|\mathcal{B}_{2}]$ is the unique $\mu$-almost everywhere  $\mathcal{B}_{2}$-measurable
function such that
\[ \int_{B}{\E[f|\mathcal{B}_{2}]d\mu}=\int_{B}{fd\mu},\quad  B \in \mathcal{B}_{2}\]
combined with  Hille's theorem \cite[Thm.~II.6]{DiestelUhl}, that for each $h \in H$ we have
\[\langle h,  \int_{B}{\E[f|\mathcal{B}_{2}]d\mu}\rangle =
\int_{B}{\langle h,\E[f|\mathcal{B}_{2}]\rangle d\mu}, \quad B \in \mathcal{B}_{2}\]
and \[\langle h, \int_{B}{f d\mu}\rangle = \int_{B}{\langle h, f \rangle d\mu}, \quad  B \in \mathcal{B}_{2}\, ,\]
 implies
 that
\[ \E[\langle h,f\rangle |\mathcal{B}_{2}]=  \langle h, \E[f|\mathcal{B}_{2}]\rangle,   \quad h \in H\]
thus implying
the following commutative diagram for all $h \in H$:
\begin{equation}
\label{fig0}
\begin{xy}
 (0,0)*+{L^{2}(\R,\mu,\mathcal{B})}="a";
(45,0)*+{L^{2}(\R,\mu,\mathcal{B}_{2})}="c";
 (0,30)*+{ L^{2}(H,\mu,\mathcal{B})}="a1";
(45,30)*+{ L^{2}(H,\mu,\mathcal{B}_{2})}="c1";
{\ar^{\E[\, \cdot\, |\mathcal{B}_{2}]} "a1";"c1"};
{\ar^{\E[\, \cdot\, |\mathcal{B}_{2}]} "a";"c"};
{\ar_{\langle h} "a1";"a"};
{\ar_{\langle h} "c1";"c"};
\end{xy}
\end{equation}

When $\mu$ is  a Gaussian measure,  the theory of Gaussian Hilbert spaces, see e.g.~Jansen \cite{janson1997gaussian}, provides
a stronger characterization of conditional expectation of the canonical random variable
$X(h):=h, h \in H$  when conditioning on a subspace and captures the full linear nature of
Gaussian conditioning. Let us assume  henceforth that $\mu$ is a centered Gaussian measure.
 Then Fernique's Theorem \cite{fernique1975regularite}, see \cite[Thm.~2.6]{da2014stochastic}, implies that the
random variable $X$ is square Bochner integrable. For any element $h \in H$, let us denote the corresponding
function $\xi_{h}:H \rightarrow \R$ defined by $\xi_{h}(h')=\langle h,h'\rangle, h' \in H$.  Then the
 the discussion above shows that for any
$h \in H$, that the real-valued random variable $\xi_{h} $ is square integrable, that is
$\xi_{h} \in  L^{2}(\R,\mu,\mathcal{B})$, for all $h \in H$.
Let
\[ \xi:H  \rightarrow  L^{2}(\R,\mu,\mathcal{B})\]
denote the resulting linear mapping
defined by
\[  h \mapsto \xi_{h} \in  L^{2}(\R,\mu,\mathcal{B}), h \in H\, .\]
It is straightforward to show that $\xi$ is injective if and only if
the covariance operator $C$ of the Gaussian measure $\mu$ is injective.
 By the definition
of a centered Gaussian vector $X$,  it follows that the law $(\xi_{h})_{*}\mu$ in $\R$ is a univariate centered Gaussian measure, that is $\xi_{h}$ is a  centered Gaussian real-valued random variable.
Consequently, let us consider the closed linear subspace
\[ H^{\mu}:=\overbar{\xi(H)} \subset  L^{2}(\R,\mu,\mathcal{B})\]
generated by  the  elements  $\xi_{h}\in  L^{2}(\R,\mu,\mathcal{B}), h \in H$.
By Jansen \cite[Thm.~I.1.3]{janson1997gaussian}, this closure
 $ H_{\mu} \subset  L^{2}(\R,\mu,\mathcal{B})$ also consists of  centered Gaussian
random variables, and since it is a closed subspace of a Hilbert space, it is a Hilbert space
and therefore  a Gaussian Hilbert space as defined in Jansen \cite{janson1997gaussian}.
Moreover, by
Jansen \cite[Thm.~8.15]{janson1997gaussian},  $H_{\mu}$ is a feature space for the Cameron-Martin
reproducing kernel Hilbert space with feature map $\xi:H \rightarrow H^{\mu}$ and reproducing kernel  the covariance operator.
For a closed Hilbert subspace, $H_{2}\subset H$, we can consider
the  closed linear subspace
\[ H_{2}^{\mu}:=\overbar{\xi(H_{2})} \subset  L^{2}(\R,\mu,\mathcal{B}_{2})\]
generated by  the  elements  $\xi_{h_{2}},h_{2} \in H_{2}$ in the same way. $ H_{2}^{\mu}$ is also a Gaussian Hilbert space and we have
the natural subspace identification
$ H_{2}^{\mu} \subset H^{\mu} $.
Since  separable Hilbert spaces are Polish,  and an
 orthonormal basis is a separating set,
it follows, see e.g.~Vakhania, Tarieladze and Chobanyan \cite[Thm.~I.1.2]{vakhania1987probability}, that for an orthonormal
basis $e_{i},i \in I$ of a separable Hilbert space, that the $\s$-algebra generated by the corresponding real-valued
functions
$\s(\{\xi_{e_{i}},i\in I\})$ is the Borel $\s$-algebra of the Hilbert space.
 Consequently, we obtain from   Janson \cite[Thm.~9.1]{janson1997gaussian}  that
for any $ h \in H$, that
\begin{eqnarray*}
\E[\xi_{h}|\mathcal{B}_{2}]&=&
\E[\xi_{h}|\s(\cup{\xi_{h_{2}}},h_{2} \in H_{2})]\\
&=& P_{H_{2}^{\mu}}\xi_{h}
\end{eqnarray*}
where $P_{H_{2}^{\mu}}:H^{\mu}\rightarrow H_{2}^{\mu}$ is orthogonal projection.
That is, if we let  $\E[\, \cdot|\mathcal{B}_{2}]: L^{2}(\R,\mu,\mathcal{B})  \rightarrow
 L^{2}(\R,\mu,\mathcal{B}_{2})
$
be the conditional expectation represented as  orthogonal projection
  and
 $\E[\, \cdot|\mathcal{B}_{2}]:  H^{\mu} \rightarrow H_{2}^{\mu}$
be the conditional expectation represented as  orthogonal projection
from the linear subspace   $H^{\mu} \subset L^{2}(\R,\mu,\mathcal{B})$
 onto the closed subspace $ H_{2}^{\mu} \subset H^{\mu} $,
we have the following commutative diagram, where $i_{H^{\mu}}:H^{\mu} \rightarrow
L^{2}(\R,\mu,\mathcal{B})$ and $i_{H_{2}^{\mu}}:H_{2}^{\mu} \rightarrow
L^{2}(\R,\mu,\mathcal{B}_{2})$  denote the closed subspace injections.
\begin{equation}
\label{fig1}
\begin{xy}
 (0,0)*+{H^{\mu}}="a";
(45,0)*+{H_{2}^{\mu}}="c";
 (0,30)*+{ L^{2}(\R,\mu,\mathcal{B})}="a1";
(45,30)*+{ L^{2}(\R,\mu,\mathcal{B}_{2})}="c1";
{\ar^{\E[\, \cdot\, |\mathcal{B}_{2}]} "a1";"c1"};
{\ar^{\E[\, \cdot\, |\mathcal{B}_{2}]} "a";"c"};
{\ar^{i_{H^{\mu}}} "a";"a1"};
{\ar^{i_{H_{2}^{\mu}}} "c";"c1"};
\end{xy}
\end{equation}
which when combined with Figure \ref{fig0}, representing the commutativity of vector projection and conditional expectation, produce the following commutative diagram for all $h \in H$:
\begin{equation}
\label{fig2}
\begin{xy}
(0,-30)*+{H}="a0";
(45,-30)*+{H_{2}}="c0";
(0,0)*+{H^{\mu}}="a";
(45,0)*+{H_{2}^{\mu}}="c";
 (0,30)*+{ L^{2}(\R,\mu,\mathcal{B})}="a1";
(45,30)*+{ L^{2}(\R,\mu,\mathcal{B}_{2})}="c1";
{\ar^{\E[\, \cdot\, |\mathcal{B}_{2}]} "a1";"c1"};
{\ar^{\E[\, \cdot\, |\mathcal{B}_{2}]} "a";"c"};
{\ar^{i_{H^{\mu}}} "a";"a1"};
{\ar^{i_{H_{2}^{\mu}}} "c";"c1"};
 (0,60)*+{ L^{2}(H,\mu,\mathcal{B})}="a2";
(45,60)*+{ L^{2}(H,\mu,\mathcal{B}_{2})}="c2";
{\ar^{\E[\, \cdot\, |\mathcal{B}_{2}]} "a2";"c2"};
{\ar_{\langle h} "a2";"a1"};
{\ar_{\langle h} "c2";"c1"};
{\ar^{\xi} "a0";"a"};
{\ar^{\xi} "c0";"c"};
\end{xy}
\end{equation}
Although  there is a natural projection map $P_{H_{2}}:H \rightarrow  H_{2}$ for the bottom of this diagram,
in general  it cannot be inserted here and maintain the commutativity of the diagram. This
 comes from the fact that there may exist an $h \in H$ such that $\xi_{h}=0$. However,
this does not imply that $\xi_{P_{H_{2}}h}=0$.

We are now prepared to obtain the main assertion.
The covariance operator
of the random variable $X$ is defined by
\begin{eqnarray*}
 \bigl\langle C h, h' \bigr\rangle&=&\E_{\mu}\bigl[\langle X,h\rangle \langle X,h'\rangle\bigr]\\
&=&\E_{\mu}\bigl[\xi_{h}\xi_{h'}\bigr], \quad  h,h'
\in H\, .
 \end{eqnarray*}
Moreover,  by the theorem of normal correlation and the commutativity of the  diagram \eqref{fig0}, the conditional  covariance operator is defined by
\begin{eqnarray*}
\bigl\langle C(X|X_{2}) h, h' \bigr\rangle&=&\E_{\mu}\Bigl[\bigl\langle X- \E[X|\mathcal{B}_{2}],h
\bigr\rangle
\bigl\langle  X-\E[X|\mathcal{B}_{2}] ,h'\bigr\rangle|\mathcal{B}_{2}\Bigr]
\\
&=&\E_{\mu}\Bigl[\bigl\langle X- \E[X|\mathcal{B}_{2}],h
\bigr\rangle
\bigl\langle  X-\E[X|\mathcal{B}_{2}] ,h'\bigr\rangle \Bigr]\\
&=&\E_{\mu}\Bigl[\bigl(\xi_{h}-\E[\xi_{h}|\mathcal{B}_{2}]\bigr)
\bigl(\xi_{h'}-\E[\xi_{h'}|\mathcal{B}_{2}]\bigr)\Bigr],\quad
h,h'
\in H\, .
\end{eqnarray*}
In terms of the Gaussian Hilbert  spaces $H^{\mu}_{2}\subset H^{\mu}$, using  the commutativity of the diagram
\eqref{fig1}
and the identification of the conditional
expectation with orthogonal projection, we conclude that
\begin{equation}
\label{id_cov}
\bigl\langle C h, h' \bigr\rangle=\langle \xi_{h},\xi_{h'}\rangle_{H^{\mu}}, \quad h, h'
\in H\,
\end{equation}
and
\begin{equation}
\label{id_condcov}
\bigl\langle C(X|X_{2}) h, h' \bigr\rangle=
\langle (I-P_{H_{2}^{\mu}}) \xi_{h},(I-P_{H_{2}^{\mu}})\xi_{h'}\rangle_{H^{\mu}}, \quad h, h'
\in H \, .
\end{equation}
Since  the orthogonal projection $P_{H_{2}^{\mu}}$ is a metric projection of $H^{\mu}$ onto $H_{2}^{\mu}$,
 we can express the dual optimization problem to the metric projection as follows:
 for any $h \in H$, using the decomposition $h=h_{1}+h_{2}$ with $h_{1}\in H_{1}, h_{2} \in H_{2}$, we decompose
$\xi_{h}=\xi_{h_{1}+h_{2}}=\xi_{h_{1}}+\xi_{h_{2}}$. Then, noting that
$(I-P_{H_{2}^{\mu}})\xi_{h_{2}}=0$, we   obtain
\begin{eqnarray*} \snorm{\xi_{h}}^{2}_{H^{\mu}}&=& \snorm{\xi_{h_{1}}+\xi_{h_{2}}}^{2}_{H^{\mu}}\\
&=& \snorm{(I-P_{H_{2}^{\mu}})(\xi_{h_{1}}+\xi_{h_{2}})}^{2}_{H^{\mu}}+ \snorm{P_{H_{2}^{\mu}}(\xi_{h_{1}}+\xi_{h_{2}})}^{2}_{H^{\mu}}\\
&=& \snorm{(I-P_{H_{2}^{\mu}})\xi_{h_{1}}}^{2}_{H^{\mu}}+ \snorm{P_{H_{2}^{\mu}}\xi_{h_{1}}+\xi_{h_{2}}}^{2}_{H^{\mu}}\, .
\end{eqnarray*}
Since in the second term on the right-hand side $P_{H_{2}^{\mu}}\xi_{h_{1}} \in H_{2}^{\nu}$, there is a sequence
$ h_{2}^{n}, n=1,...$ such that the corresponding sequence $\xi_{h_{2}^{n}}$ converges
to $-P_{H_{2}^{\mu}}\xi_{h_{1}}$ in $ L^{2}(\R,\mu,\mathcal{B})$ and therefore $H^{\mu}$,
we conclude that
\[ \snorm{(I-P_{H_{2}^{\mu}})\xi_{h_{1}}}^{2}_{H^{\mu}}=\inf_{h_{2}\in H_{2}}{\snorm{\xi_{h_{1}}+\xi_{h_{2}}}^{2}_{H^{\mu}}}\,  \, .\]
From the identifications \eqref{id_cov} and  \eqref{id_condcov}, we conclude that
\[\bigl\langle C(X|X_{2}) h_{1}, h_{1} \bigr\rangle =\inf_{h_{2}\in H_{2}}{\bigl\langle C(X)(h_{1}+h_{2}),
h_{1}+ h_{2} \bigr\rangle},\, . \]
Therefore,   Anderson and Trapp \cite[Thm.~6]{anderson1975shorted} implies the assertion
\[ C(X|X_{2})=\mathcal{H}_{1}\bigl(C\bigr)\, .\]
The assertion in the non-centered case follows by simple translation.

\subsection{Proof of Theorem \ref{thm_short}}
Since  $P_{H_{2}}A^{n}P_{H_{2}} =P_{H_{2}^{n}}A^{n}P_{H_{2}^{n}}$, the range of $P_{H_{2}}A^{n}P_{H_{2}}$ is finite dimensional,  and therefore closed, so that it follows from
Corach, Maestripieri and Stojanoff \cite[Lem.~3.8]{corach2006projections} that
$A^{n}$ is compatible with $H_{2}$ for all $n$.

Now we utilize the approximation  results  of Butler and Morley \cite{butler1988note} for
the shorted operator.
By \cite[Lem.~1]{butler1988note}, for $c \in H$  and for fixed $n$, it follows that there exists
a sequence $y_{m}^{n}\in H_{2}, m=1,\ldots$ and  a real number $M$ such that
\begin{equation*}
\begin{array}{cccc}
A^{n}_{11}c+A^{n}_{12}y^{n}_{m} & \rightarrow &\mathcal{H}_{1}(A^{n})c\, , &  \qquad   m \rightarrow \infty\\
A^{n}_{21}c+A^{n}_{22}y^{n}_{m} & \rightarrow &  0\, ,  &   \qquad  m \rightarrow \infty\\
\langle A^{n}_{22}y^{n}_{m}, y^{n}_{m} \rangle & \leq &  M\, ,  & \qquad  \forall m \, .
\end{array}
\end{equation*}
Since $A^{n}_{11}=A_{11}$, $A^{n}_{12}= A_{12}P_{H_{2}^{n}}$, $A^{n}_{21}= P_{H_{2}^{n}}A_{21}$, and
$A^{n}_{22}= P_{H_{2}^{n}}A_{22}P_{H_{2}^{n}}$ this can be written as
\begin{equation*}
\begin{array}{cccc}
A_{11}c+A_{12}P_{H_{2}^{n}}y^{n}_{m} & \rightarrow & \mathcal{H}_{1}(A^{n})c \, , & \qquad   m \rightarrow \infty\\
P_{H_{2}^{n}}A_{21}c+P_{H_{2}^{n}}A_{22}P_{H_{2}^{n}}y^{n}_{m} & \rightarrow &  0\, , &   \qquad  m \rightarrow \infty\\
\langle A_{22}P_{H_{2}^{n}}y^{n}_{m}, P_{H_{2}^{n}}y^{n}_{m} \rangle & \leq &  M\, , &\qquad  \forall m\, .
\end{array}
\end{equation*}
Since  these equations only depend on $P_{H_{2}^{n}}y^{n}_{m}$ we can further assume
that $P_{(H_{2}^{n})^{\perp}}y^{n}_{m}=0, \, m =1,\ldots,$ where
 $P_{(H_{2}^{n})^{\perp}}$ is the orthogonal projection
onto $(H_{2}^{n})^{\perp} \subset H_{2}$. That is, we can assume that
$P_{H_{2}^{n}}y^{n}_{m}=y^{n}_{m}, m=1,\ldots$ and therefore
\begin{equation}
\label{eq_zzzz}
\begin{array}{cccc}
A_{11}c+A_{12}y^{n}_{m} & \rightarrow & \mathcal{H}_{1}(A^{n})c \, , & \qquad   m \rightarrow \infty\\
P_{H_{2}^{n}}A_{21}c+P_{H_{2}^{n}}A_{22}y^{n}_{m} & \rightarrow &  0\, , &   \qquad  m \rightarrow \infty\\
\langle A_{22}y^{n}_{m}, y^{n}_{m} \rangle & \leq &  M\, , &\qquad  \forall m\, .
\end{array}
\end{equation}
It follows from
$\mathcal{H}_{1}(A^{n}) \leq A^{n}$
that $\snorm{\sqrt{\mathcal{H}_{1}(A^{n})}}\leq \snorm{\sqrt{A^{n}}}$
 for the unique square root, guaranteed to exist by
 Riesz and Sz.-Nagy \cite[Sec.~104]{riesz1955functional}. Consequently,
 Conway \cite[Prop.~II.2.7]{conway1990course} implies that
 $\snorm{\mathcal{H}_{1}(A^{n})}\leq \snorm{A^{n}}$  for all $n$ and since
$\snorm{A^{n}}\leq \snorm{A}$ for all $n$ it follows that
$\snorm{\mathcal{H}_{1}(A^{n})}\leq \snorm{A}$  for all $n$. Consequently, the sequence
 $\mathcal{H}_{1}(A^{n})c$ is bounded. Therefore there exists a weakly convergent subsequence.
Let $n'$ denote the index of any weakly convergent subsequence, so that
\begin{equation}
\label{eq_bbbb}
\mathcal{H}_{1}(A^{n'})c \xrightarrow{\omega} d',\quad n' \rightarrow \infty \,
\end{equation}
for some $d'$ depending on the subsequence.
Now the strong convergence of the lefthand side to the righthand side in \eqref{eq_zzzz} is maintained for the subsequence
$n'$ and, since for the subsequence the  first term on the righthand side converges weakly to $d'$,
it follows that we can define a monotonically increasing function $m(n')$ and use it to define a new sequence
$\hat{y}^{n'}:=y^{n'}_{m(n')}$  such
that
\begin{equation}
\label{eq_zzzz2}
\begin{array}{cccc}
A_{11}c+A_{12}\hat{y}^{n'} & \xrightarrow{\omega} & d'\, , & \qquad   n' \rightarrow \infty\\
P_{H_{2}^{n'}}A_{21}c+P_{H_{2}^{n'}}A_{22}\hat{y}^{n'}  & \rightarrow &  0\, , &   \qquad  n' \rightarrow \infty\\
\langle A_{22}\hat{y}^{n'} , \hat{y}^{n'}  \rangle & \leq &  M\, , &\qquad  \forall n'\, .
\end{array}
\end{equation}
Since $P_{H_{2}^{n}}$ is strongly convergent to $P_{H_{2}}$ it follows that
 $P_{H_{2}^{n'}}$ is strongly convergent to $P_{H_{2}}$, so that
$P_{H_{2}^{n'}}A_{21}c$ converges to $A_{21}c$ and $P_{H_{2}^{n'}}A_{22}\hat{y}^{n'}$ converges to $-A_{21}c$. Moreover,
by Reid's inequality
\cite[Cor.~2]{reid1951symmetrizable}
we have
\begin{equation}
\label{ineq_reid}
  \snorm{A_{22}\hat{y}^{n'}}^{2}_{H_{2}}\leq \snorm{A_{22}}\langle A_{22}\hat{y}^{n'}, \hat{y}^{n'}\rangle
\leq \snorm{A_{22}} M\, ,
\end{equation}
for all $n'$, so that the sequence $A_{22}\hat{y}^{n'}$ is bounded.
Since weak convergence of a bounded sequence on a separable Hilbert space is equivalent to
the convergence with respect to each element of any orthonormal basis, it
  follows that
$A_{22}\hat{y}^{n'}$ is weakly convergent to $-A_{21}c$. From \eqref{eq_zzzz2},
we obtain
\begin{equation}
\label{eq_abc}
\begin{array}{cccc}
A_{12}\hat{y}^{n'} & \xrightarrow{\omega} & d'-A_{11}c\, , & \qquad   n' \rightarrow \infty\\
A_{22}\hat{y}^{n'}& \xrightarrow{\omega} & -A_{21}c \, , &  \qquad  n' \rightarrow \infty\, .
\end{array}
\end{equation}

 From Kakutani's \cite{kakutani1938weak} generalization of the Banach-Saks Theorem
it follows that we can select a subsequence $\grave{n}$ of $n'$ such that
 the Cesaro means of $A_{22}\hat{y}^{\grave{n}}$ and $A_{12}\hat{y}^{\grave{n}}  $  converge strongly  in \eqref{eq_abc}.
That is, if we consider  the Cesaro means
\[z^{\grave{n}}=\frac{1}{\grave{n}}\sum_{i=1}^{\grave{n}}{\hat{y}^{\grave{n}}}\]
we have
\begin{equation*}
\begin{array}{cccc}
A_{12}z^{n'} & \rightarrow & d'-A_{11}c\, , & \qquad   n' \rightarrow \infty\\
A_{22}z^{n'}& \rightarrow & -A_{21}c \, , &  \qquad  n' \rightarrow \infty \, .
\end{array}
\end{equation*}
Since $A_{22}\geq 0$ it follows that the function $y \mapsto \langle A_{22}y,y\rangle$ is convex, so  that $\langle A_{22} z^{\grave{n}}, z^{\grave{n}} \rangle \leq M $  for all $\grave{n}$, so that
\begin{equation*}
\label{eq_zzzz3}
\begin{array}{cccc}
A_{11}c+A_{12}z^{\grave{n}} & \rightarrow & d'\, , & \qquad   \grave{n} \rightarrow \infty\\
A_{21}c+A_{22}z^{\grave{n}}& \rightarrow & 0 \, , &  \qquad  \grave{n} \rightarrow \infty\\
\langle A_{22}z^{\grave{n}}, z^{\grave{n}} \rangle & \leq & M\, , & \qquad  \forall \grave{n}\, .
\end{array}
\end{equation*}
It therefore follows from the from the main result of
Butler and Morley \cite[Thm.~1]{butler1988note} that
\[d'  =\mathcal{H}_{1}(A)c\, .\]
Consequently, by \eqref{eq_bbbb}, we obtain  that
\begin{equation}
\label{eq_bbbb2}
\mathcal{H}_{1}(A^{n'})c \xrightarrow{\omega} \mathcal{H}_{1}(A)c,\quad n' \rightarrow \infty \, .
\end{equation}
Since this limit  is independent of the chosen weakly converging subsequence, it follows,
 see e.g.~Zeidler \cite[Prop.~10.13]{zeidler1989nonlinear}, that the full sequence weakly converges
to the same limit, that is we have
 \begin{equation}
\label{eq_bbbb3}
\mathcal{H}_{1}(A^{n})c \xrightarrow{\omega} \mathcal{H}_{1}(A)c,\quad n \rightarrow \infty \, ,
\end{equation}
and since $c$ was arbitrary we conclude that
\[ \mathcal{H}_{1}(A)= \omega\mbox{-}\\lim_{n\rightarrow\infty}{ \mathcal{H}_{1}(A^{n})}\, .\]

\subsection{Proof of Theorem \ref{thm_main}}
Let us first establish the assertion when $C$ is compatible with $H_{2}$.
 Consider the operator $\hat{C}:H \rightarrow H$ defined by
\[  \hat{C}:=\mathcal{H}_{1}(C) +P_{2} CP_{2}\, .\]
Since $C$ is compatible with $H_{2}$, there exists an oblique projection
 $Q \in \mathcal{P}(C,H_{2})$, and
Corach, Maestripieri and Stojanoff~\cite[Prop.~4.2]{corach2001schur}
 asserts that for
  $E:=1-Q$,  we have
\begin{equation}
\label{e1}
\mathcal{H}_{1}(C)=CE=E^{*}CE\, .
\end{equation}
Since  $Q^{*}C=CQ$ it follows that $E^{*}C=CE$, and
since $Q$ is a projection, it follows that $QE=EQ=0$ and that $E$ is a projection.
Moreover, since $R(Q)=H_{2}$ it follows that $ker(E)=H_{2}$, so that we obtain
$P_{2}Q=Q$ and $EP_{1}=E$ and therefore $Q^{*}P_{2}=Q^{*}$ and $P_{1}E^{*}=E^{*}$.
Consequently,  we obtain
\begin{eqnarray*}
 (P_{1}+Q)^{*}\hat{C}(P_{1}+Q)&=& (P_{1}+Q)^{*}(E^{*}CE+P_{2}CP_{2})(P_{1}+Q)\\
&=& (P_{1}+Q)^{*}(E^{*}CE+P_{2}CQ)\\
&=&E^{*}CE+ Q^{*}CQ \\
&=& CE+CQ \\
&=& C\, ,
\end{eqnarray*}
that is,
\begin{equation}
\label{id_covariance}
C=(P_{1}+Q)^{*}\hat{C}(P_{1}+Q)\, .
\end{equation}

Since $Q$ is a projection onto $H_{2}$, it follows that $P_{1}+Q$ is lower triangular in its  partitioned
representation and therefore the fundamental pivot produces  an explicit, and most importantly continuous, inverse.
Indeed, if we use the
partition representation
\[Q= \begin{pmatrix}
  0& 0\\
\hat{Q}& 1
 \end{pmatrix}\, ,
\]
we see that
\[ (P_{1}+Q)= \begin{pmatrix}
  1 & 0  \\
   \hat{Q} & 1
 \end{pmatrix}
\]
from which we conclude that
\[ (P_{1}+Q)^{-1}= \begin{pmatrix}
  1 &  0  \\
   -\hat{Q} & 1
 \end{pmatrix}
\]
Without partitioning, using
 $P_{1}Q=0$ and $QP_{2}=P_{2}$,   we obtain
\begin{eqnarray*}
(2-P_{1}-Q)(P_{1}+Q) &=&2P_{1}+2Q -(P_{1}^{2}+P_{1}Q+QP_{1}+Q^{2})\\
&=&2P_{1} +2Q- P_{1} -P_{1}Q-QP_{1}-Q\\
&=&P_{1}+Q-QP_{1}\\
&=&P_{1}+QP_{2}\\
&=&P_{1} +P_{2}\\
&=& 1
\end{eqnarray*}
and so confirm that
 \begin{equation}
\label{eq_inv}
(P_{1}+Q)^{-1}=2-P_{1}-Q\, .
\end{equation}

Following the proof of Hairer, Stuart, Voss, and Wiber
\cite[Lem.~4.3]{hairer2005analysis},
let $\mathcal{N}(m,C)$ denote the Gaussian measure with mean $m$ and covariance operator $C$ and consider the transformation
\[(P_{1}+Q)^{-*}:H \rightarrow H\, ,\]
where we use the notation $A^{-*}$ for $(A^{-1})^{*}=(A^{*})^{-1}$.  From \eqref{id_covariance} we obtain
\begin{equation}
\label{covariance2}
(P_{1}+Q)^{-*}C(P_{1}+Q)^{-1}=\hat{C}
\end{equation}
 so that the transformation  law
for Gaussian measures, see
Maniglia and Rhandi \cite[Ch.~1,~Lem.~1.2.7]{maniglia2004gaussian}, implies  that
\[\bigl((P_{1}+Q)^{-*}\bigr)_{*}\mathcal{N}(m,C)=\mathcal{N}\bigl((P_{1}+Q)^{-*}m,\hat{C}\bigr)\, . \]

Since
\[ (P_{1}+Q)^{-1}= \begin{pmatrix}
  1 &  0  \\
   -\hat{Q} & 1
 \end{pmatrix}
\]
we obtain
\[ (P_{1}+Q)^{-*}
= \begin{pmatrix}
  1 &  -\hat{Q}^{*}  \\
   0 & 1
 \end{pmatrix}
\]
and therefore
\[ (P_{1}+Q)^{-*}m=
\begin{pmatrix}
  m_{1}-\hat{Q}^{*}m_{2}  \\
  m_{2}
 \end{pmatrix}\, .
\]
 Since the partition representation of $\hat{C}$ is
\[\hat{C}=  \begin{pmatrix}
(\mathcal{H}_{1}(C))_{11}  &  0  \\
  0 &  C_{22}
 \end{pmatrix}
\]
 the components of the corresponding Gaussian random variable are uncorrelated and therefore independent.
That is, we have
\[\mathcal{N}\bigl((P_{1}+Q)^{-*}m,\hat{C}\bigr)=
\mathcal{N}\bigl(m_{1}-\hat{Q}^{*}m_{2},(\mathcal{H}_{1}(C))_{11}\bigr)
\mathcal{N}\bigl(m_{2},C_{22}\bigr)
\, .\]
This independence facilitates the computation of the conditional measure as follows. Let
$X=(X_{1},X_{2})$ denote the random variable associated with the Gaussian measure
$\mathcal{N}(m,C)$ and consider the transformed random variable
$Y=(P_{1}+Q)^{-*}X$ with the product law
$\mathcal{N}\bigl(m_{1}-\hat{Q}^{*}m_{2},(\mathcal{H}_{1}(C))_{11}\bigr)\mathcal{N}\bigl(m_{2},C_{22}\bigr)$.
Then,
\begin{eqnarray*}
Y_{1}&=& X_{1}-\hat{Q}^{*}X_{2}\\
Y_{2}&=&X_{2}
\end{eqnarray*} can be used to compute
 the conditional expectation  as
\begin{eqnarray*}
\E[X_{1}|X_{2}]&=&\E[X_{1}-\hat{Q}^{*}X_{2}|X_{2}]+\E[\hat{Q}^{*}X_{2}|X_{2}]\\
&=&\E[Y_{1}|Y_{2})]+\E[\hat{Q}^{*}X_{2}|X_{2}]\\
&=&\E[Y_{1}]+\hat{Q}^{*}X_{2}\, ,
\end{eqnarray*}
obtaining
\begin{equation}
\E[X_{1}|X_{2}]=
\E[Y_{1}]+\hat{Q}^{*}X_{2}\, ,
\end{equation}
so that we conclude that
\[\E[X_{1}|X_{2}]=m_{1}+\hat{Q}^{*}(X_{2}-m_{2})\, .\]
A similar calculation obtains the covariance
\begin{equation}
\label{id_cov}
C(X|X_{2})=\mathcal{H}_{1}(C),
\end{equation}
 thus establishing the assertion in the compatible case.

For the general case, we do not assume that
 $C$ is compatible with $H_{2}$.   Consider
 an ordered orthonormal basis
for $H_{2}$, let $H^{n}_{2}$ denote the span of the first $n$ basis elements,
 let $P^{n}:=P_{H_{1}}+P_{H^{n}_{2}}$ denote the orthogonal projection onto $H_{1}\oplus H^{n}_{2}$, and
consider the sequence of Gaussian measures $\mu_{n}:=P^{n}_{*}\mu$ with the  mean $P^{n}m$ and
  covariance operators
\[C^{n}:=P^{n}CP^{n},\quad n=1,\ldots\]
As asserted in Theorem \ref{thm_short}, $C^{n}$  is
 compatible with $H_{2}$ for all $n$,
 and the  sequence
$\mathcal{H}_{1}(C^{n})$   converges weakly to $\mathcal{H}_{1}(C)$.
Let $C(X_{1}|X_{2}^{n})$ and $C(X_{1}|X_{2})$ denote the conditional covariance operators associated
with the measure $\mu$.
Then we will show that
$C(X_{1}|X_{2}^{n})=\mathcal{H}_{1}(C^{n})$, so that
the assertion regarding the conditional covariance operators
is established if we demonstrate that the sequence of conditional covariance operators $C(X_{1}|X_{2}^{n})$
 converges weakly to $C(X_{1}|X_{2})$.

To both ends,
consider the Lebesgue-Bochner space $L^{2}(H,\mu,\mathcal{B})$  space of (equivalence classes) of $H$-valued
Borel measurable functions on $H$ whose squared norm
\[ \snorm{f}^{2}_{L^{2}(H,\mu,\mathcal{B})}:=\int_{H}{\snorm{f(x)}^{2}_{H}d\mu(x)}\ \]
is integrable.
Since
 Fernique's Theorem \cite{fernique1975regularite}, see \cite[Thm.~2.6]{da2014stochastic}, implies that the
random variable $X$ is square Bochner integrable, it follows that the Gaussian random variables
$P^{n}X$ are also square Bochner integrable with respect to $\mu$.
Let us denote
$\mathcal{B}_{2}:= \{ (H_{1},T):  T \in \mathcal{B}(H_{2})\}$ and
$\mathcal{B}^{n}_{2}:=\{ (H_{1},T^{n},(H_{2}^{n})^{\perp}):  T^{n} \in \mathcal{B}(H^{n}_{2})\}$, and
let $\mu_{n}:=P^{n}_{*}\mu$ denote the image under the projection. $\mu_{n}$ is a Gaussian measure
on $H$
with mean $P^{n}m$ and covariance $C^{n}$.

Now consider a function $f:H \rightarrow H$ which is  Bochner square integrable with respect to
$\mu$ and satisfies $f\circ P^{n}=f$. Then, using the change of variables formula for
Bochner integrals, see Bashirov \cite[Thm.~2.26]{bashirov2003partially}, along with the fact that
  $(P^{n})^{-1}\mathcal{B}_{2}=\mathcal{B}^{n}_{2}$, and using the fact that for an arbitrary
$\mathcal{B}^{n}_{2}$-measurable function $g$ we have $g=g\circ P^{n}$, it follows that
 for $ A \in \mathcal{B}_{2}$,
 we have
\begin{eqnarray*}
 \int_{A}{fd\mu_{n}}&=& \int_{(P^{n})^{-1}A}{f\circ P^{n}d\mu}\\
&=& \int_{(P^{n})^{-1}A}{fd\mu}\\
&=&\int_{(P^{n})^{-1}A}{\E_{\mu}\bigl[f|(P^{n})^{-1}\mathcal{B}_{2}\bigr]d\mu}\\
&=&\int_{(P^{n})^{-1}A}{\E_{\mu}\bigl[f|\mathcal{B}^{n}_{2}\bigr]d\mu}\\
&=&\int_{(P^{n})^{-1}A}{\E_{\mu}\bigl[f|\mathcal{B}^{n}_{2}\bigr]\circ P^{n}d\mu}\\
&=&\int_{A}{\E_{\mu}\bigl[f|\mathcal{B}^{n}_{2}\bigr]d\mu_{n}}\\
\end{eqnarray*}
we obtain
\begin{equation}
\label{eq_martingale}
 \E_{\mu_{n}}\bigl[f|\mathcal{B}_{2}\bigr]=  \E_{\mu}\bigl[f|\mathcal{B}^{n}_{2}\bigr]\, ,
\end{equation}
and conclude that the sequence $ \E_{\mu_{n}}\bigl[f|\mathcal{B}_{2}\bigr], n=1 \ldots $ is a
 martingale corresponding to the increasing
family of $\s$-algebras $\mathcal{B}^{n}_{2}$.
Moreover, it is easy to see that \eqref{eq_martingale} holds for
real valued functions $f:H \rightarrow R$ which are square integrable with respect to
$\mu$ and satisfy $f\circ P^{n}=f$.
With the choice $f:=X_{1}$, we clearly have $X_{1}\circ P^{n}=X_{1}$, so that if we denote
$X^{n}_{2}:=P^{n}X_{2}$,  we conclude that the sequence
\begin{equation}
\label{id_martingale}
\E_{\mu_{n}}[X_{1}|X_{2}]=\E_{\mu}[X_{1}|X^{n}_{2}]\, ,\quad  n=1,\ldots
\end{equation}
is a martingale.
Since conditional expectation is a contraction, it follows that
 the $L^{2}$ norm of all the conditional
expectations are uniformly bounded by the $L^{2}$ norm of  $X$.
Then by the Martingale Convergence Theorem of  Diestel and Uhl \cite[Cor.~V.2.2] {DiestelUhl},
$\E_{\mu_{n}}\bigl[X_{1}|X_{2}\bigr]$ converges to
$\E_{\mu}\bigl[X_{1}|X_{2}\bigr]$ in $L^{2}(H,\mu,\mathcal{B})$.

For the conditional covariance operators,
observe that \eqref{id_martingale} implies that
\begin{equation}
\label{id_martingale2}
X-\E_{\mu_{n}}[X|X_{2}]= X_{1}-\E_{\mu}[X_{1}|X^{n}_{2}]\,
\end{equation}
for all $n$, so that
 for $h_{1}, h_{2} \in H$,  we have
\begin{eqnarray*}
\langle C_{\mu_{n}}(X|X_{2})h_{1},h_{2}\rangle
&:=&\E_{\mu_{n}}\Bigl[\bigl\langle X-\E_{\mu_{n}}[X|X_{2}], h_{1}\bigr\rangle\bigl\langle X-\E_{\mu_{n}}[X|X_{2}] ,  h_{2}\bigr\rangle\big|X_{2}\Bigr]\\
&=&\E_{\mu_{n}}\Bigl[\bigl\langle X_{1}-\E_{\mu}[X_{1}|X^{n}_{2}], h_{1}\bigr\rangle\bigl\langle X_{1}
-\E_{\mu}[X_{1}|X^{n}_{2}] ,
h_{2}\bigr\rangle\big|X_{2}\Bigr]
\end{eqnarray*}
and since the integrand  $ f:=\bigl\langle X_{1}-\E_{\mu}[X_{1}|X^{n}_{2}], h_{1}\bigr\rangle\bigl\langle X_{1}
-\E_{\mu}[X_{1}|X^{n}_{2}] ,
h_{2}\bigr\rangle$ satisfies $f\circ P^{n}=f$, it follows from \eqref{eq_martingale}
that
\begin{eqnarray*}
&&\E_{\mu_{n}}\Bigl[\bigl\langle X_{1}-\E_{\mu}[X_{1}|X^{n}_{2}], h_{1}\bigr\rangle\bigl\langle X_{1}
-\E_{\mu}[X_{1}|X^{n}_{2}] ,
h_{2}\bigr\rangle\big|X_{2}\Bigr]\\
&=&\E_{\mu}\Bigl[\bigl\langle X_{1}-\E_{\mu}[X_{1}|X^{n}_{2}], h_{1}\bigr\rangle\bigl\langle X_{1}
-\E_{\mu}[X_{1}|X^{n}_{2}] ,
h_{2}\bigr\rangle\big|X^{n}_{2}\Bigr]
\end{eqnarray*}
so that using the theorem of normal correlation,
we obtain
\begin{equation}
\label{id_cov33}
\langle C_{\mu_{n}}(X|X_{2})h_{1},h_{2}\rangle=\E_{\mu}\Bigl[\bigl\langle X_{1}-\E_{\mu}[X_{1}|X^{n}_{2}], h_{1}\bigr\rangle\bigl\langle X_{1}
-\E_{\mu}[X_{1}|X^{n}_{2}] ,
h_{2}\bigr\rangle\Bigr]\, .
\end{equation}
Since the theorem of normal correlation also shows that
\begin{eqnarray*}
\langle C_{\mu}(X|X_{2})h_{1},h_{2}\rangle&:=&
\E_{\mu}\Bigl[\bigl\langle X-\E_{\mu}[X|X_{2}], h_{1}\bigr\rangle\bigl\langle X-\E_{\mu}[X|X_{2}] ,  h_{2}\bigr\rangle \big|X_{2}\Bigr]\\
&=&\E_{\mu}\Bigl[\bigl\langle X-\E_{\mu}[X|X_{2}], h_{1}\bigr\rangle\bigl\langle X-\E_{\mu}[X|X_{2}] ,  h_{2}\bigr\rangle \Bigr]\\
&=&
\E_{\mu}\Bigl[\bigl\langle X_{1}-\E_{\mu}[X_{1}|X_{2}], h_{1}\bigr\rangle\bigl\langle X_{1}-\E_{\mu}[X_{1}|X_{2}] ,  h_{2}\bigr\rangle \Bigr]\, ,
\end{eqnarray*}
 the difference in the covariances can be decomposed as
\begin{eqnarray*}
&&\langle C_{\mu_{n}}(X_{1}|X_{2})h_{1},h_{2}\rangle-\langle C_{\mu}(X_{1}|X_{2})h_{1},h_{2}\rangle \\
&=& \E_{\mu}\Bigl[\bigl\langle X_{1}-\E_{\mu}[X_{1}|X^{n}_{2}], h_{1}\bigr\rangle\bigl\langle X_{1}-\E_{\mu}[X_{1}|X^{n}_{2}] ,  h_{2}\bigr\rangle \Bigr]\\
&&-\E_{\mu}\Bigl[\bigl\langle X_{1}-\E_{\mu}[X_{1}|X_{2}], h_{1}\bigr\rangle\bigl\langle X_{1}-\E_{\mu}[X_{1}|X_{2}] ,  h_{2}\bigr\rangle \Bigr]\\
&=& \E_{\mu}\Bigl[\bigl\langle \E_{\mu}[X_{1}|X_{2}]- \E_{\mu}[X_{1}|X^{n}_{2}], h_{1}\bigr\rangle\bigl\langle X_{1},  h_{2}\bigr\rangle \Bigr]+
\E_{\mu}\Bigl[\bigl\langle X_{1},h_{1}\bigr\rangle\bigl\langle \E_{\mu}[X_{1}|X_{2}]- \E_{\mu}[X_{1}|X^{n}_{2}] ,  h_{2}\bigr\rangle \Bigr]\\
&&+\E_{\mu}\Bigl[\bigl\langle \E_{\mu}[X_{1}|X^{n}_{2}], h_{1}\bigr\rangle\bigl\langle \E_{\mu}[X_{1}|X^{n}_{2}] ,  h_{2}\bigr\rangle \Bigr]-\E_{\mu}\Bigl[\bigl\langle \E_{\mu}[X_{1}|X_{2}], h_{1}\bigr\rangle\bigl\langle \E_{\mu}[X_{1}|X_{2}] ,  h_{2}\bigr\rangle \Bigr]
\end{eqnarray*}
where the last term can be decomposed as
\begin{eqnarray*}
&&\E_{\mu}\Bigl[\bigl\langle \E_{\mu}[X_{1}|X^{n}_{2}], h_{1}\bigr\rangle\bigl\langle \E_{\mu}[X_{1}|X^{n}_{2}] ,  h_{2}\bigr\rangle \Bigr]-\E_{\mu}\Bigl[\bigl\langle \E_{\mu}[X_{1}|X_{2}], h_{1}\bigr\rangle\bigl\langle \E_{\mu}[X_{1}|X_{2}] ,  h_{2}\bigr\rangle \Bigr]\\
=&& \E_{\mu}\Bigl[\bigl\langle \E_{\mu}[X_{1}|X^{n}_{2}]-\E_{\mu}[X_{1}|X_{2}], h_{1}\bigr\rangle\bigl\langle \E_{\mu}[X_{1}|X^{n}_{2}] ,  h_{2}\bigr\rangle \Bigr]\\
 &+& \E_{\mu}\Bigl[\bigl\langle \E_{\mu}[X_{1}|X_{2}], h_{1}\bigr\rangle\bigl\langle \E_{\mu}[X_{1}|X^{n}_{2}]-\E_{\mu}[X_{1}|X_{2}] ,  h_{2}\bigr\rangle \Bigr]\, .
\end{eqnarray*}
Then since conditional expectation is a contraction on $L_{2}(H,\mu,\mathcal{B})$ it follows
that
\newline
$\snorm{\E_{\mu}[X_{1}|X_{2}]}_{L_{2}(H,\mu,\mathcal{B})} \leq \snorm{X_{1}}_{L_{2}(H,\mu,\mathcal{B})}$
and
$\snorm{\E_{\mu}[X_{1}|X^{n}_{2}]}_{L_{2}(H,\mu,\mathcal{B})} \leq \snorm{X_{1}}_{L_{2}(H,\mu,\mathcal{B})}$
for all $n$. Moreover, since $\E_{\mu}\bigl[X_{1}|X^{n}_{2}\bigr]$ converges to
$\E_{\mu}\bigl[X_{1}|X_{2}\bigr]$ in $L^{2}(H,\mu,\mathcal{B})$ it follows,
 see e.g.~\cite[Lem.~11.45]{AliprantisBorder:2006}, that
$\langle \E_{\mu}\bigl[X_{1}|X^{n}_{2}\bigr],h\rangle$ converges to
$\langle \E_{\mu}\bigl[X_{1}|X_{2}\bigr], h\rangle$   in $L^{2}(\R,\mu,\mathcal{B})$ for all
$h \in H$. Therefore, the Cauchy-Schwartz inequality applied four times in the above decomposition implies
that
\[ \lim_{n\rightarrow \infty}\langle C_{\mu_{n}}(X|X_{2})h_{1},h_{2}\rangle
=\langle C_{\mu}(X_{1}|X_{2})h_{1},h_{2}\rangle, \quad h_{1}, h_{2} \in H \]
so that   we obtain
\[ C_{\mu}(X|X_{2})= \omega\mbox{-}\\lim_{n\rightarrow\infty}{C_{\mu_{n}}(X|X_{2})}\, .\]
Since $C^{n}$ is compatible with $H_{2}$ for all $n$, and the compatible case demonstrated  in
\eqref{id_cov} that
\begin{equation}
\label{id_cov2}
C_{\mu_{n}}(X|X_{2})=\mathcal{H}_{1}(C^{n})
\end{equation}
for all $n$, and
   Theorem \ref{thm_short} asserts that
\[\mathcal{H}_{1}(C)= \omega\mbox{-}\\lim_{n\rightarrow\infty}{\mathcal{H}_{1}(C^{n})}\, ,\]
we conclude that
$C_{\mu}(X|X_{2})=\mathcal{H}_{1}(C)$, establishing the assertion regarding the covariance operators.

For the means, observe that
since $\mu$  is a probability measure,
it follows that $X$ and therefore $X_{1}$ lie in the Lebesgue-Bochner space
 $L^{1}(H,\mu,\mathcal{B})$, and since  by  Diestel and Uhl \cite[Thm.~V.1.4]{DiestelUhl} the conditional expectation operators
are also contractions on $L^{1}(H,\mu,\mathcal{B})$ it also follows that
$\E_{\mu_{n}}\bigl[X_{1}|X_{2}\bigr]$ converges to
$\E_{\mu}\bigl[X_{1}|X_{2}\bigr]$ in $L^{1}(H,\mu,\mathcal{B})$.
Therefore,  Diestel and Uhl \cite[Thm.~V.2.8]{DiestelUhl}  implies that
$\E_{\mu_{n}}\bigl[X_{1}|X_{2}\bigr]$ converges to
$E_{\mu}\bigl[X_{1}|X_{2}\bigr]$ a.e.-$\mu$. Let
the conditional means $\E_{\mu}\bigl[X|X_{2}\bigr]$ be
denoted by $\E_{\mu}\bigl[X|X_{2}\bigr]=m_{t}, t\in H_{2}$.  Then, since
\[P^{n}m =\begin{pmatrix}
  m_{1} \\
P_{H^{n}_{2}}m_{2}
 \end{pmatrix}\, ,
\] is
the mean of the measure $\mu_{n}$,
 the assertion in the
 compatible case demonstrated that the conditional means
  $\E_{\mu_{n}}\bigl[X|X_{2}]=m^{n}_{t}, t\in H_{2}$
are
\[m^{n}_{t}= \begin{pmatrix}
  m_{1}+\hat{Q}_{n}^{*}(t-P_{H^{n}_{2}}m_{2}) \\
t
 \end{pmatrix}\, .
\]
 Since the conditional means $\E_{\mu_{n}}\bigl[X_{1}|X_{2}\bigr]$ converge to the conditional means
$\E_{\mu}\bigl[X_{1}|X_{2}\bigr]$ a.e.-$\mu$ amounts to
$m^{n}_{t} \rightarrow m_{t}$ for $\mu$-almost every $t$,
 the first assertion regarding the means is also proved. Now suppose
that $Q_{n}$ eventually becomes the special element $Q_{n}=Q_{C^{n},H_{2}}$ defined
near \eqref{id_special}. Then, by definition, $R(\hat{Q}_{n})\subset \overbar{R(C^{n}_{22})}$ so that
$ker(\hat{Q}^{*}_{n})\supset R(C^{n}_{22})^{\perp}$, but since
$C^{n}_{22}=\Pi_{2}C^{n}\Pi_{2}^{*}=\Pi_{2}P^{n}CP^{n}\Pi_{2}^{*}=\Pi_{2}P_{H^{n}_{2}}CP_{H^{n}_{2}}\Pi_{2}^{*}$,
 it follows that $R(C^{n}_{22})\subset H_{2}^{n}$
and therefore $R(C^{n}_{22})^{\perp}\supset (H^{n}_{2})^{\perp}$ so that
$ker(\hat{Q}^{*}_{n}) \supset (H^{n}_{2})^{\perp}$. Therefore
$\hat{Q}^{*}_{n}P_{H^{n}_{2}}=\hat{Q}^{*}_{n}$, so that the final assertion follows from the previous.

\subsection{Proof of Corollary \ref{cor_short}}
By Mourier's Theorem, there exists a Gaussian measure $\mu$ on $H$ with mean $0$ and
covariance operator $C:=A$.
Looking at the end of the proof of Theorem \ref{thm_main},
 since conditional expectation is a contraction on $L_{2}(H,\mu,\mathcal{B})$ it follows
that
$\snorm{\E_{\mu}[X_{1}|X_{2}]}_{L_{2}(H,\mu,\mathcal{B})} \leq \snorm{X_{1}}_{L_{2}(H,\mu,\mathcal{B})}$
and
$\snorm{\E_{\mu}[X_{1}|X^{n}_{2}]}_{L_{2}(H,\mu,\mathcal{B})} \leq \snorm{X_{1}}_{L_{2}(H,\mu,\mathcal{B})}$
for all $n$. Therefore, for $h \in H$,  it follows from the Cauchy-Schwartz inequality that
$\snorm{\langle \E_{\mu}[X_{1}|X^{n}_{2}], h\rangle}_{L_{2}(\mathbb{R},\mu,\mathcal{B})} \leq \snorm{X_{1}}_{L_{2}(H,\mu,\mathcal{B})}$ and $\snorm{\langle \E_{\mu}[X_{1}|X_{2}], h\rangle}_{L_{2}(\mathbb{R},\mu,\mathcal{B})} \leq \snorm{X_{1}}_{L_{2}(H,\mu,\mathcal{B})}$ for all $n$, uniformly for $h \in H$ with $\snorm{h}_{H}\leq 1$.
Therefore, the Cauchy-Schwartz inequality applied four times in the  decomposition at the end of the
proof of Theorem \ref{thm_main} implies
that
\[ \lim_{n\rightarrow \infty}\langle C_{\mu_{n}}(X|X_{2})h_{1},h_{2}\rangle
=\langle C_{\mu}(X_{1}|X_{2})h_{1},h_{2}\rangle, \quad h_{1}, h_{2} \in H \]
uniformly for $h_{1},h_{2} \in H$ with $\snorm{h_{1}}_{H}\leq 1$ and $\snorm{h_{2}}_{H}\leq 1$. Therefore, it
follows from Halmos \cite[Prob.~107]{halmos1982hilbert} that the sequence of covariance operators converges
\[ C_{\mu_{n}}(X|X_{2}) \rightarrow
C_{\mu}(X|X_{2})\]
in the uniform operator topology.

According to Maniglia and Rhandi \cite[Ch.~1,~Lem.~1.1.4]{maniglia2004gaussian} or
Da Prato and Zabczyk \cite[Prop.~2.16]{da2014stochastic}, for a Gaussian measure $\mu$ with mean
$0$ and covariance operator $C$, we have
\begin{equation*}
\label{id_trace}tr(C)=\E_{\mu}{\snorm{X}^{2}}\, .
\end{equation*}
From \eqref{id_cov33},  by
 shifting to the center, we obtain  that
\begin{equation*}
\label{id_cov33a}
tr\bigl(C_{\mu_{n}}(X|X_{2})\bigr)=\E_{\mu}\Bigl[\snorm{X_{1}-\E_{\mu}[X_{1}|X^{n}_{2}]}^{2}
\Bigr]\,
\end{equation*}
and
\begin{eqnarray*}
tr\bigl(C_{\mu}(X|X_{2})\bigr)=
\E_{\mu}\Bigl[\snorm{X_{1}-\E_{\mu}[X_{1}|X_{2}]}^{2} \Bigr]\, ,
\end{eqnarray*}
and
therefore the difference is
\begin{eqnarray*}
&&tr\bigl(C_{\mu_{n}}(X|X_{2})\bigr)-tr\bigl(C_{\mu}(X|X_{2})\bigr) \\
&=& \E_{\mu}\Bigl[\snorm{X_{1}-\E_{\mu}[X_{1}|X^{n}_{2}]}^{2}
\Bigr]-E_{\mu}\Bigl[\snorm{X_{1}-\E_{\mu}[X_{1}|X_{2}]}^{2} \Bigr]\\
&=& \E_{\mu}\Bigl[\bigl\langle \E_{\mu}[X_{1}|X^{n}_{2}] -\E_{\mu}[X_{1}|X_{2}],\,
\E_{\mu}[X_{1}|X^{n}_{2}] +\E_{\mu}[X_{1}|X_{2}]-2X_{1}\bigr\rangle \Bigr]\,.
\end{eqnarray*}
Therefore, the Cauchy-Schwartz inequality,  the $L^{2}$ convergence of
$\E_{\mu}[X_{1}|X^{n}_{2}]$ to $ \E_{\mu}[X_{1}|X_{2}]$, and the uniform $L^{2}$ boundedness of
$\E_{\mu}[X_{1}|X^{n}_{2}]$, $ \E_{\mu}[X_{1}|X_{2}]$ and $X_{1}$,
implies that
\[ \lim_{n\rightarrow \infty}{tr\bigl(C_{\mu_{n}}(X|X_{2})\bigr)}=tr\bigl(C_{\mu}(X|X_{2})\bigr)\, .\]
Since $C_{\mu_{n}}(X|X_{2}) \rightarrow C_{\mu}(X|X_{2})$ in the uniform operator topology,
it follows from
Kubrusly \cite{Kubrusly}, see \cite[Thm.~2 \& Sec.~4]{kubrusly1985convergence},
that
$C_{\mu_{n}}(X|X_{2}) \rightarrow C_{\mu}(X|X_{2})$ in the trace norm topology.
Since  \eqref{id_cov2} asserts that
$C_{\mu_{n}}(X|X_{2})=\mathcal{H}_{1}(C^{n})$ and Theorem \ref{thm_main}
asserts that $C_{\mu}(X|X_{2})=\mathcal{H}_{1}(C)$, the identification  $A:=C$ completes the proof.
\section*{Acknowledgments}
The authors gratefully acknowledge this work supported by the Air Force Office of Scientific Research under Award Number
FA9550-12-1-0389 (Scientific Computation of Optimal Statistical Estimators).
\addcontentsline{toc}{section}{Acknowledgments}

\newpage

\addcontentsline{toc}{section}{References}
\bibliographystyle{plain}
\bibliography{refs}

\end{document}